\documentclass[a4paper,12pt]{article}
\usepackage{amsmath, amssymb, amsthm, latexsym, amsopn, amsfonts}
\usepackage{amscd,stmaryrd}
\usepackage[cp1251]{inputenc}
\usepackage{mathtext}
\usepackage[english]{babel}
\usepackage{indentfirst}

\begin{document}

\begin{center}
{\bf\large Multiplicative Automorphisms of Incidence Algebras}
\end{center}
\begin{center}
Evgenii Kaigorodov\footnote{Gorno-Altaisk State University, e-mail: gazetaintegral@gmail.com},
Piotr Krylov\footnote{National Research Tomsk
State University, e-mail: krylov@math.tsu.ru},
Askar Tuganbaev\footnote{National Research University MPEI; Lomonosov Moscow State University; e-mail: tuganbaev@gmail.com}
\end{center}

\sloppy

\textbf{Abstract. }Let $I(X,R)$ be the incidence algebra of the preordered set $X$ over the ring $R$. In the case of a finite connected partially ordered set $X$, we prove that the subgroup of inner multiplicative automorphisms is a direct factor of the group of multiplicative automorphisms of the algebra $I(X,R)$. As a consequence, we obtain several matching criteria of the subgroup of inner multiplicative automorphisms with the group of multiplicative automorphisms. 


\textbf{MSC2020 database 16W20, 16S50}

The study is supported by Russian Scientific Foundation. 







\section{Introduction}\label{section1}

The incidence algebras of preordered (in particular, partially ordered) sets represent an important and very characteristic example of rings of functions. In this paper, we study multiplicative automorphisms of the incidence algebra $I(X,R)$ of a preordered set $X$ over the ring $R$. Automorphisms of incidence algebras have long attracted the attention of specialists. For information on this subject, you can refer to the book \cite{SpiO97} and the papers \cite{BruFS15}, \cite{BruL11}, \cite{DroK07}, \cite{Khr10}. Multiplicative automorphisms are interesting in themselves. They play an essential role in the study of arbitrary automorphisms of the algebra $I(X,R)$; moreover, they act as one of the types of its standard automorphisms (see \cite{BruL11}, \cite{DroK07}, \cite{SpiO97}). In addition, multiplicative automorphisms are used to describe involutions of the algebra $I(X,R)$ (see \cite{BruFS12}, \cite{BruFS14}, \cite{BruL11}).

Incidence algebras were defined by Rota in his familiar paper \cite{Rot64} as a tool for solving some problems of combinatorics and, above all, for studying generalizations of the M\"obius inversion formula in number theory in a unified way.
Over time, incidence algebras themselves, regardless of their applications in combinatorics and other areas of mathematics, have turned out to be a meaningful algebraic object. Many works are dedicated to them, including the book \cite{SpiO97}.

Section 2 of the given paper contains some definitions and auxiliary  results. In Section 3, we first give several general facts on multiplicative automorphisms. The main result of the section and the whole paper is Theorem 3.6 which states that if $X$ is a finite connected partially ordered set, then the subgroup of inner multiplicative automorphisms is a direct factor in the group of multiplicative automorphisms of the algebra $I(X,R)$. 
As a consequence of this theorem, several conditions have been found under which this multiplicative automorphism is inner.

We only consider associative rings with non-zero identity element. For some ring (or an algebra) $S$, $U(S)$ is its group of invertible elements, $C(S)$ is its center, $M(n,\,S)$ is the ring of all square $n$-order matrices with values in $S$, $M(m\!\times\!n,\,S)$ is the group of all rectangular matrices of order $m\!\times\!n$ with values in $S$. Further, $\mathrm{Aut}\,S$ is the automorphism group of the algebra $S$, $\mathrm{In}(\mathrm{Aut}\,S)$ is a subgroup of its inner automorphisms. The semidirect product of two groups $G$ and $H$ is denoted by $G\leftthreetimes H$.

\section{Preliminaries}\label{section2}

We briefly outline some initial information about preordered sets (one can get acquainted with them in more detail in the book \cite{SpiO97}).

Let $X$ be an arbitrary set and let $\leqslant$ be a reflexive transitive relation on $X$. In such a case, the system $\langle X,\leqslant \rangle$ is called a \textsf{preordered set} and $\leqslant$ is a \textsf{preorder} on $X$. If the relation $\leqslant$ is also antisymmetric, then $\langle X,\leqslant\rangle$ is a \textsf{partially ordered set}.

Further, we assume that $\langle X,\leqslant\rangle$ is a preordered set. For any two elements $x,\,y\in X$, we denote by $[x,\,y]$ the subset $\{z\in X\;|\;x\leqslant z\leqslant y\}$. It is called an \textsf{interval} in $X$. An interval of the form $[x,x]$ is denoted by $[x]$. We note two useful properties of intervals.

\textbf{(a)} For any $y,\,z\in [x]$, we have the equality $[y,\,z]=[x]$;

\textbf{(b)} If $x< y$, then $s<t$ for arbitrary elements $s\in [x]$ and $t\in [y]$.

We define a binary relation $\sim$ on $X$ by setting $x\sim y$ $\Leftrightarrow$ $x\leqslant y$ and $y\leqslant x$. It is clear that $\sim$ is an equivalence relation on $X$. The corresponding equivalence classes are of the form $[x]$ for every $x\in X$. It follows from the property $(b)$ that the pre-order relation $\leqslant$ is agreed with the equivalence relation $\sim$. Consequently, the induced  relation $\leqslant$ appears on the factor set $\overline{X}=X/\!\!\sim$ and $\langle \overline{X},\leqslant \rangle$ is a partially ordered set.

With the preordered set $X$ and the partially ordered set $\overline{X}$, we can associate a directed graph (we do not take into account the loops that arise in this case). It is more convenient to proceed from the set $\overline{X}$. When necessary, we consider $\overline{X}$ as a simple graph related to a directed graph $\overline{X}$. In this case, we use standard notions of graph theory such that a \textsf{semi-path} and its \textsf{length}, \textsf{spanning tree}, \textsf{cycle} and  \textsf{cyclomatic number}.

We agree that all intervals in $X$ are finite. In such a case, $X$ is called a \textsf{locally finite} preordered set.

The \textsf{length} of the interval in a locally finite partially ordered set is called the greatest of lengths of chains in this interval.

We agree in the future to denote by $x$ elements of the partially ordered set $\overline{X}$, i.e., the equivalence classes of the form $[x]$. In other words, we use some representative of the class $[x]$ to represent $[x]$. This arrangement is correct and should not lead to confusion. In a specific situation, it is always clear which elements of the set ($\overline{X}$ or $X$) we are talking about.

Further, let the letter $R$ denote an algebra over some commutative ring $T$. However, the ring $T$ itself is hardly used explicitly.

An incidence algebra is some function ring. Let $\langle X,\leqslant\rangle$ be an arbitrary locally finite preordered set. We set $I(X,R)=\{f\colon X\!\!\times\!\!X\to R\;|\; f(x,\,y)=0,\, \text{ if }\, x\not\leqslant y\}$. Functions are added pointwise. The product of functions $f$ and $g$ from $I(X,R)$ is given by the formula
$$
(fg)(x,\,y)=\sum_{x\leqslant z\leqslant y}f(x,\,z)\cdot g(z,\,y)\eqno(1)
$$
for any $x,\,y\in X$. Since $X$ is a locally finite set, we can write $z\in X$ instead of $x\leqslant z\leqslant y$ in $(1)$. For any $t\in T$ and $x,\,y\in X$, we also assume $(tf)(x,\,y)=tf(x,\,y)$. As a result, we obtain a $T$-algebra $I(X,R)$ called the \textsf{incidence algebra} or the \textsf{ring incidence} of the preordered set $X$ over the ring $R$. In what follows, we denote by $A$ the specific algebra $I(X,R)$.

We define some special functions in $I(X,R)$. For $x\in X$, we set $e_{[x]}(t,t)=1$ for all $t\in [x]$ and $e_{[x]}(z,\,y)=0$ for remaining pairs $(z,\,y)$. The system $\{e_{[x]}\;|\;x\in X\}$ consists of pairwise orthogonal idempotents which are central in $L$ (the ring $L$ is defined in the next paragraph). In accordance with the convention of denoting the class $[x]$ by some of its representatives, we will write $e_x$ instead of $e_{[x]}$.

We define a subring $L$ and an ideal $M$ in $A$. We set $L=\{f\in A\;|\;f(x,\,y)=0,\, \text{ if }\, x\not\sim y\}$ and $M=\{f\in A\;|\;f(x,\,y)=0,\, \text{ if }\, x\sim y\}$. We have the direct sum $A=L\oplus M$ of $T$-modules, i.e., the ring $A$ is a splitting extension of the ideal $M$ with the use of the subring $L$. The ideal $M$ is a natural $L$-$L$- bimodule. In addition, $M$ is a non-unital algebra.

The ideal $M$ is contained in the Jacobson radical of the algebra $A$. Consequently, the element $1+d$ invertible in $A$ for every $d\in M$ \cite[Theorem 1.2.3]{SpiO97}.

Let we have an arbitrary interval $[x]$. We denote by $R_{[x]}$ the set of functions $f\in A$ such that $f(z,\,y)=0$ if $z\not\sim x$ or $y\not\sim x$. Similar to the case of idempotents $e_x$, we write $R_x$ instead of $R_{[x]}$. There are the equalities $R_x=e_xAe_x=e_xLe_x.$ We conclude that $R_x$ is a ring with identity element $e_x$. If we pass to restrictions of functions in $A$ to $[x]\!\times\![x]$, then, in fact, $R_x$ is an algebra of all functions $[x]\!\times\![x]\to R$ with pointwise addition and product of convolution type as in $(1)$. We take some numbering of the interval $[x]\colon [x]=\{x_1,\ldots,x_n\}$. After this, if the function $f\in R_x$ is associated with the matrix $(f(x_i,x_j))$, then we obtain an algebra isomorphism $R_x\cong M(n,R)$. Now we take two distinct intervals $[x]$, $[y]$ and set
$$
M_{xy}=\{f\in A\;|\;f(s,t)=0,\, \text{ if }\, s\not\sim x\, \text{ or } \, t\not\sim y\}.
$$
Then $M_{xy}=e_xAe_y$ and, therefore, $M_{xy}$ is an $R_x$-$R_y$-bimodule. In addition, the equality $M_{xy}=e_xMe_y$ is true.

We clarify that $M_{xy}=0$ for $x\not\leqslant y$. For $x< y$, there exists the canonical isomorphism
$$
M_{xy}\cong M(n\!\times\!m,R),\quad n=\big|[x]\big|,\; m=\big|[y]\big|
$$
with respect to above isomorphisms $R_x\cong M(n,R)$, $R_y\cong M(m,R)$. After identifications of all algebras $R_x$ with $M(n,R)$ and bimodules $M_{xy}$ with $M(n\!\times\!m,R)$, it becomes clear that ring actions $R_x$ and $R_y$ on $M_{xy}$ are ordinary multiplications of matrices. It is clear also that $M_{xy}$ is an $L$-$L$-bimodule. The action $L$ on $M_{xy}$ reduces to the action $R_x$ on the left and $R_y$ on the right.

We note once again that by the indices in $R_x$ and $M_{xy}$ we mean $[x]$ and $[x][y]$, respectively (see above).

The product $\!\!\!\prod\limits_{x,\,y\in X}\!\!M_{xy}$ is 
an $L$-$L$-bimodule. Namely, if $f\in L$, $(g_{xy})\in\!\!\!\prod\limits_{x,\,y\in X}\!\!M_{xy}$, then $f(g_{xy})=(f_xg_{xy})\, \text{ and }\,(g_{xy})f=(g_{xy}f_y)$
where $f_x=e_xfe_x$ and $f_y=e_yfe_y$.

In the bimodule $\!\!\!\prod\limits_{x,\,y\in X}\!\!M_{xy}$, we define a multiplication through the formula
$(g_{xy})\cdot (h_{xy})=(d_{xy}),$
where $d_{xy}=\!\!\!\sum\limits_{x\leqslant z\leqslant y}\!\!g_{xz}\cdot h_{zy}$. After which this bimodule becomes a (non-unital) algebra.

\textbf{Proposition 2.1.} There exist canonical algebra isomorphisms and also $L$-$L$-bimodule isomorphisms $L\cong\!\prod\limits_{x\in X}R_{x}$, and also an algebra isomorphism $M\cong\!\!\!\prod\limits_{x,\,y\in X}\!\!M_{xy}$.

\textbf{Proof.} We define a mapping $\omega\colon L\to\!\!\prod\limits_{x\in X}\!R_x$ by setting $\omega(f)=(f_x)$ for every $f\in L$, where $f_x=e_xfe_x$. Then $\omega$ is an algebra isomorphism.

The mapping $\varepsilon\colon M\to\!\!\!\prod\limits_{x,\,y\in X}\!\!M_{xy}$, $\varepsilon(g)=(g_{xy})$, where $g_{xy}=e_xge_y$, is an $L$-$L$-bimodule isomorphism and an algebra isomorphism.~$\square$

In what follows, we don't distinguish objects corresponding to each other under isomorphisms $\omega$ and $\varepsilon$.

\textbf{Remark 2.2.} We can consider an algebra $A$ as a function  algebra and as an abstract ring presented as a splitting extension $L\oplus M$. These two approaches can be called a "functional approach" and an "abstract approach"; of course, they are equivalent. We use the both of these approaches.

Every invertible element $v$ of the algebra $A$ provides an automorphism $\mu_v$ of $A$ with the use of the following rule: $\mu_v(a)=v^{-1}av$, $a\in A$. Such an automorphism is called an inner automorphism defined by the element $v$. All inner automorphisms form a normal subgroup $\mathrm{In}(\mathrm{Aut}\,A)$ of the group $\mathrm{Aut}\,A$. We denote by $\mathrm{In_1}(\mathrm{Aut}\,A)$ (correspondingly, $\mathrm{In_0}(\mathrm{Aut}\,A)$) the subgroup of inner automorphisms of the algebra $A$ defined by invertible elements of the form $1+d$, $d\in M$ (correspondingly, by invertible elements of the algebra $L$). The first subgroup is normal in $\mathrm{Aut}\,A$. By the use of the semidirect decomposition $U(A)=(1+M)\leftthreetimes U(L)$, we can verify that the decomposition $\mathrm{In}(\mathrm{Aut}\,A)=\mathrm{In_1}(\mathrm{Aut}\,A)\leftthreetimes\mathrm{In_0}(\mathrm{Aut}\,A)$ exists (also see \cite[Section 4]{KryT21}).

\section{Multiplicative Automorphisms}\label{section3}

In this section, we consider the ring $L$ and the bimodules $M$, $M_{xy}$ defined above. It follows from Remark 2.2 that $A=I(X,R)=L\oplus M$. In fact, we do not impose any conditions on the preordered set $X$, the ring $R$, and the algebra $A$ which is equal to $I(X,R)$.

Let $\varphi$ be an arbitrary automorphism of the algebra $A$. 
Based on the $T$-module direct sum $A=L\oplus M$, the automorphism $\varphi$ (as well as any $T$-endomorphism) can be associated in the usual way with a $2$-order square matrix $\begin{pmatrix}\alpha& \gamma\\ \delta &\beta\end{pmatrix}$. We will only deal with matrices of the form $\begin{pmatrix}1&0\\ 0&\beta\end{pmatrix}$. If the automorphism $\varphi$ corresponds to such a matrix, then $\beta$ is an automorphism of the (non-unital) algebra $M$ and an automorphism of the $L$-$L$-bimodule $M$. Conversely, let some mapping $\beta\colon M\to M$ satisfy the above properties. We obtain an automorphism $\varphi$ of the algebra $A$ by setting $\varphi(a+c)=\begin{pmatrix}1&0\\ 0&\beta\end{pmatrix}\begin{pmatrix}a\\c\end{pmatrix}=a+\beta(c)$ for any $a\in L$, $c\in M$. All these assertions are directly verified. Their validity also follows from more general considerations in Section 3 of of the paper \cite{KryT21}. In what follows, we not distinguish automorphisms of the algebra $A$ and the matrices corresponding to them.

Automorphisms of the algebra $A$ of the form $\begin{pmatrix}1&0\\ 0&\beta\end{pmatrix}$ are said to be \textsf{multiplicative}. All multiplicative automorphisms form a normal subgroup in $\mathrm{Aut}\,A$. We denote it by $\mathrm{Mult}\,A$.

Inner automorphisms of the algebra $A$, corresponding to invertible central elements of the ring $L$, are said to be \textsf{fractional}. We denote by $\mathrm{Mult_0}\,A$ the subgroup of all fractional automorphisms. In the book \cite{SpiO97}, the definitions of multiplicative and fractional automorphisms are seemingly different definitions from those given here (there is a text on this topic at the end of the section). The following result establishes a connection between groups $\mathrm{Mult}\,A$ and $\mathrm{Mult_0}\,A$.

\textbf{Proposition 3.1.} There are equalities 
$$
\mathrm{Mult}\,A\cap\text{In}(\text{Aut}\,A)=\mathrm{Mult}\,A\cap\text{In}_0(\text{Aut}\,A)=\mathrm{Mult_0}\,A.
$$

\textbf{Proof.} We verify a non-obvious inclusion $\mathrm{Mult}\,A\cap\text{In(Aut}\,A)\subseteq\mathrm{Mult_0}\,A$. We take some automorphism $\mu=\mu_1\mu_0$, where $\mu\in\mathrm{Mult}\,A$, $\mu_1\in\text{In}_1(\text{Aut}\,A)$, $\mu_0\in\text{In}_0(\text{Aut}\,A)$. Let $\mu_1$ be defined by the element $1+d$, $d\in M$, and let $\mu_0$ be defined by the element $v\in U(L)$. Then $\mu$ is defined by the element $v+vd$, and $\mu^{-1}$ is defined by the element $v^{-1}+d'$, where $d'\in M$.

For any element $a\in L$, we have
$$
\mu(a)=v^{-1}av+(v^{-1}avd+d'av+d'avd).
$$
Since $\mu\in\mathrm{Mult}\,A$, we have that $v^{-1}av=a$ and $av=va$. Therefore, $v\in C(L)$. The expression in parentheses is zero. In particular, $v^{-1}e_xvd+d^\prime e_xv+d^\prime e_xvd=0$, $e_xv^{-1}e_xvd+e_xd^\prime e_xv+e_xd^\prime e_xvd=e_xd=0$ for every $x\in X$. Therefore, $d=0$, whence we have $\mu_1=1$ and $\mu=\mu_0$. Since $v$ is an invertible central element, $\mu\in\mathrm{Mult_0}\,A$.~$\square$

It can be concluded that fractional automorphisms coincide with inner multiplicative automorphisms (see also Proposition 3.4).

In this section, we consider various information about groups $\mathrm{Mult}\,A$ and $\mathrm{Mult_0}\,A$. If $X$ is a finite connected set, then the subgroup $\mathrm{Mult_0}\,A$ can be calculated exactly (Corollary 3.7) and it is a direct factor of the group $\mathrm{Mult}\,A$ (Theorem 3.6).

The following fact is standard and essentially known (see also \cite[Proposition 13.2]{KryT21}).

\textbf{Proposition 3.2.} Let $P=M(n,R)$, $Q=M(m,R)$ and $V=M(n\!\times\!m,R)$. Endomorphisms (correspondingly, automorphisms) of the $P$-$Q$-bimodule $V$ coincide with multiplications by central elements (correspondingly, invertible central elements) of the ring $R$.

Let we have a multiplicative automorphism $\begin{pmatrix}1&0\\ 0&\beta\end{pmatrix}$. For every element $c\in M_{xy}$, we have the equalities $\beta(c)=\beta(e_xce_y)=e_x\beta(c)e_y$. Consequently, $\beta$ leaves in place every bimodule $M_{xy}$.
We denote by $\beta_{xy}$ the restriction $\beta|_{M_{xy}}$. Since $\beta$ is an automorphism of the algebra $M$, the equality
$$
\beta_{xy}(ab)=\beta_{xz}(a)\beta_{zy}(b)\eqno (2)
$$
is true for any $x,\,z,\,y\in X$ with $x<z<y$ and all $a\in M_{xz}$, $b\in M_{zy}$ . Due to the fact that $\beta$ is an automorphism of the $L$-$L$-bimodule $M$ at the same time, we have that $\beta_{xy}$ is an automorphism of the $R_x$-$R_y$-bimodule $M_{xy}$ for any $x,\,y$ with $x<y$.

Further, it follows from Proposition 3.2 that there exists an element $c_{xy}\in C(U(R))$ with property $\beta_{xy}(g)=c_{xy}g$, $g\in M_{xy}$. Taking into account the equality $M_{xz}M_{zy}=M_{xy}$, we obtain from $(2)$ the equality
$$
c_{xy}=c_{xz}c_{zy},\eqno (3)
$$
where $x<z<y$. Thus, the automorphism $\psi\in\mathrm{Mult}\,A$ can be associated with the system of invertible central elements $c_{xy}$ ($x,\,y\in X$, $x<y$) of the ring $R$. Conversely, every system of elements
$\{c_{xy}\in C(U(R))\;|\;x<y\}$, satisfying to equality $(3)$, provides a multiplicative automorphism $\psi$. Namely, for the element $g=(g_{xy})\in M$, we set $\psi(g)=(c_{xy}g_{xy})$. In this case we have to take into account the rule of multiplication of elements in $M$ (see Section 2).

We formulate everything just stated above as follows.

\textbf{Proposition 3.3.}

\textbf{1.} There is a one-to-one correspondence between multiplicative automorphisms and systems of elements $\{c_{xy}\in C(U(R))\;|\;x<y\}$ such that $c_{xy}=c_{xz}c_{zy}$ for all $x,\,z,\,y\in X$ with $x<z<y$.

\textbf{2.} The group $\mathrm{Mult}\,A$ is embedded in the product of $\mathcal{M}$ copies of the group $C(U(R))$, where $\mathcal{M}=\big|\{(x,\,y)\;|\;x<y\}\big|$.

A more formal point of view on the situation under consideration is as follows. For every pair $(x,\,y)$ of elements in $X$ with $x<y$, we take a copy $U_{xy}$ of the group $C(U(R))$. Then the group $\mathrm{Mult}\,A$ is isomorphic to a subgroup of the product $\prod\limits_{x<y}\!U_{xy}$, consisting of vectors $(u_{xy})$, the elements of which satisfy the equality $u_{xy}=u_{xz}u_{zy}$ if $x<z<y$.

What can be said about fractional automorphisms, i.e. on automorphisms from the subgroup $\mathrm{Mult_0}\,A$? Let a fractional automorphism $\psi=\begin{pmatrix}1&0\\ 0&\beta\end{pmatrix}$ be defined by a invertible central element $v=(v_x)\in L$, where $v_x\in C(U(R))$. For any $x$, $y$ with $x<y$ and every $g\in M_{xy}$, we have the equalities
$$
\beta(g)=v^{-1}gv=v_x^{-1}v_yg.
$$
Consequently, the system of elements $\{c_{xy}\;|\;x<y\}$ from Proposition 3.3, corresponding to the automorphism $\psi$, consists of all elements of the form $v_x^{-1}v_y$, where $x<y$.

We have the following result.

\textbf{Proposition 3.4.} Let a multiplicative automorphism $\psi$ correspond to a system of elements 
$c_{xy}$ $(x,\,y\in X,\; x<y)$ satisfying equality $(2)$. The following assertions are equivalent.

\textbf{1)} $\psi$ is a fractional automorphism.

\textbf{2)} $\psi$ is an inner automorphism.

\textbf{3)} There exists an element $v=(v_x)\in L$ with $v_x\in C(U(R))$ such that is true the equality $c_{xy}=v_x^{-1}v_y$ for all $x,\,y$ $(x<y)$.

\textbf{Everywhere else we assume that $X$ is a finite connected set. For technical convenience, we also assume that $X$ is partially ordered.}

Note that the incidence algebra $I(X,R)$ for a finite set $X$ is sometimes called the \textsf{ring of structural matrices}. You can get acquainted with such rings in the papers \cite{Coe14}, \cite{SmiW94}. Structural matrix rings are one type of a large class of formal matrix rings. The book \cite{KryT17} is dedicated to the latter.

Let we have an arbitrary set $\{c_{xy}\;|\;x,\,y\in X,\; x<y\}$ of central invertible elements of the ring $R$. Based on this set, we assign weights to the edges and the half-paths of the graph $X$ (see \cite{BruFS15}). If $(x,\,y)$ is an edge, then we set $w(x,\,y)=c_{xy}$ for $x<y$ and $w(x,\,y)=c_{yx}^{-1}$ for $x>y$. For a semi-path $P=z_1z_2\ldots z_{k+1}$ in $X$, we set
$$
w(P)=\prod_{i=1}^kw(z_i,z_{i+1}).
$$
If $\psi$ is some multiplicative automorphism, then we obtain weights of edges and semi-paths of the graph $X$. They are induced by the system $\{c_{xy}\;|\;x<y\}$ corresponding to the automorphism $\psi$ in the sense of Proposition 3.3. Let
$\psi$ be a fractional automorphism (i.e., inner multiplicative automorphism; see Proposition 3.4) which is defined by an invertible  element $v=(v_x)$. It is easy to verify that for any two elements $x,\,y\in X$, the weight of any semi-path from $x$ to $y$ is equal to $v_x^{-1}v_y$ (thus, it does not depend on specific semi-path).

We fix some spanning tree $T$ of the graph $X$. For every edge $(x,\,y)$ in $T$ with $x<y$, we take an element $c_{xy}\in C(U(R))$. 
Similar to what was done above, we use the elements $c_{xy}$ to assign the weight $w$ to edges and semi-paths in $T$. Let $(x,\,y)$ be an edge in $X$ which is not contained in $T$ and $x<y$. We set $c_{xy}=w(P)$, where $P$ is a (unique) semi-path in $T$ from $x$ to $y$. As a result, we obtain a system of elements $\{c_{xy}\in C(U(R))\;|\;x,\,y\in X,\, x<y\}$ and weights of edges and semi-paths in $X$.

\textbf{Lemma 3.5.} The constructed system of elements $\{c_{xy}\;|\; x<y\}$ satisfies equalities $(3)$. Consequently, it induces a multiplicative automorphism $\psi$ in the sense of Proposition 3.3.

\textbf{Proof.} Let $x,\,z,\,y$ be elements with $x<z<y$. We denote by $P_{xz}$, $P_{zy}$ and $P_{xy}$ semi-paths in $T$ from $x$ to $z$, from $z$ to $y$ and from $x$ to $y$, correspondingly. Since semi-pathes in $T$ are unique, we have $P_{xy}=P_{xz}P_{zy}$. Therefore, $w(P_{xy})=w(P_{xz})\cdot w(P_{zy})$ and $c_{xy}=c_{xz}c_{zy}$. It follows from Proposition 3.3 that the system $\{c_{xy}\;|\; x<y\}$ corresponds to a quite definite automorphism $\psi\in\mathrm{Mult}\,A$.~$\square$

Since $X$ is a partially ordered set, the bimodules $M_{xy}$ are copies of the ring $R$ in our case. We can define analogues of matrix units. Let $x,\,y\in X$ and $x<y$. We denote by $e_{xy}$ a function $X\times X\to R$ such that $e_{xy}(s,t)=1$ if $s=x$, $t=y$ and $e_{xy}(s,t)=0$ for all remaining pairs $(s,t)$. The functions $e_{xy}$ satisfy the following property: if $x<z<y$, then $e_{xz}e_{zy}=e_{xy}$.

Note that we have chosen the spanning tree $T$ of the graph $X$. We denote by $\mathrm{Mult_1}\,A$ the subgroup in $\mathrm{Mult}\,A$ consisting of automorphisms $\psi$ such that $\psi(e_{xy})=e_{xy}$ for all edges $(x,\,y)\in T$ with $x<y$. Equivalently, $c_{xy}=1$ for the specified edges $(x,\,y)$. Here $\{c_{xy}\in C(U(R))\;|\; x<y\}$ is the system from Proposition 3.3. We also note that this proposition implies that the group $\mathrm{Mult}\,A$ is commutative.

\textbf{Theorem 3.6.} There exists a group direct decomposition $\mathrm{Mult}\,A=\mathrm{Mult_1}\,A\!\times\!\mathrm{Mult_0}\,A$.

\textbf{Proof.} Let $\psi\in\mathrm{Mult}\,A$ and let $\{c_{xy}\in C(U(R))\;|\; x<y\}$ be the corresponding system of elements from Proposition 3.3. For every vertex in $X$, we now define some central invertible element.

We take some root $x_0$ in the tree $T$. Let $y$ be some vertex in $T$ adjacent with $x_0$. If $x_0<y$, then we set $v_{x_0}=1$ and $v_y=v_{x_0}c_{x_0y}$. If $x_0>y$, then we set $v_y=v_{x_0}c_{yx_0}^{-1}$. Then we take some vertex $z$, which is adjacent with $y$ in $T$, and do the same. Namely, let $v_z=v_yc_{yz}$ for $y<z$ and let $v_z=v_{y}c_{zy}^{-1}$ for $y>z$. We continue in this manner until we reach some hanging vertex in $T$. We obtain a semi-path in $T$ from $x_0$ to this hanging vertex. If there are vertices left in $T$ that do not belong to this semi-path, then let $s$ be the first vertex in the semi-path with degree exceeding two. Starting from the vertex $s$, we similarly form another half-path, and so on. As a result, each vertex $x$ in $T$ has an element $v_x\in C(U(R))$ associated with it.

We take an invertible central element $v=(v_x)$ in $L$. Let $\psi_0$ be the inner automorphism induced by $v$. Then $\psi_0\in\mathrm{Mult_0}\,A$, since $v\in C(U(L))$. We set $\psi_1=\psi\psi_0^{-1}$. If $(x,\,y)\in T$ and $x<y$, then $\psi(e_{xy})= c_{xy}e_{xy}$, by Proposition 3.2, and $\psi_0(e_{xy})=v_x^{-1}v_ye_{xy}$. By construction of elements $v_x$, the equality $c_{xy}=v_x^{-1}v_y$ is true. Consequently, $\psi$ and $\psi_0$ act the same on elements $e_{xy}$ with $(x,\,y)\in T$. Therefore, the inclusion $\psi_1\in\mathrm{Mult_1}\,A$ is true. Thus, $\psi=\psi_1\psi_0$, where $\psi_1\in\mathrm{Mult_1}\,A$, $\psi_0\in\mathrm{Mult_0}\,A$.

It remains to verify that $\mathrm{Mult_1}\,A\cap \mathrm{Mult_0}\,A=\langle 1\rangle$. Let $\xi\in \mathrm{Mult_1}\,A\cap\mathrm{Mult_0}\,A$. Then $\xi(e_{xy})=e_{xy}$ for every edge $(x,\,y)\in T$, where $x<y$. In the same time, the automorphism $\xi$ is defined as an inner automorphism by some central invertible element $v=(v_x)$ of the ring $L$. Therefore, we have the equality $\xi(e_{xy})=v_x^{-1}v_ye_{xy}$. Consequently, $v_x^{-1}v_y=1$ and $v_x=v_y$. Now let an edge $(x,\,y)$ be not contained $T$. There exists a semi-path in $T$ from $x$ to $y$. By the use of this semi-path, we obtain the equality $v_x=v_y$ from only the above. Therefore, we obtain the equality $\xi=1$. Therefore, $\mathrm{Mult_1}\,A\cap\mathrm{Mult_0}\,A=\langle 1 \rangle$ which completes the proof.~$\square$

For a graph $X$, we denote by $m(X)$ and $\lambda(X)$ the number of edges and the cyclomatic number of $X$, respectively.

\textbf{Corollary 3.7.} There is a group isomorphism
$$
\mathrm{Mult_0}\,A\cong\!\!\!\!\!\!\!\prod_{m(X)-\lambda(X)}\!\!\!\!\!\!\!C(U(R)).
$$

\textbf{Proof.} We denote by the letter $G$ the product in the right part of the isomorphism. We define a mapping $g\colon\mathrm{Mult}\,A\to G$ as follows. Let $\psi\in\mathrm{Mult}\,A$ and let $\{c_{xy}\;|\; x<y\}$ be the system of elements from Proposition 3.3 corresponding to $\psi$. The mapping $g$ associates the automorphism $\psi$ with the vector $(c_{xy})$, where $(x,\,y)\in T$ and $x<y$. It is clear that $g$ is a group homomorphism. Surjectivity of $g$ follows from Lemma 3.5. It is also clear that $\text{Ker}\,g=\mathrm{Mult_1}\,A$. We have an isomorphism $\mathrm{Mult}\,A/\mathrm{Mult_1}\,A\cong G$. On the other hand, it follows from Theorem 3.6 that there is an isomorphism $\mathrm{Mult}\,A/\mathrm{Mult_1}\,A\cong \mathrm{Mult_0}\,A$.~$\square$

\textbf{Corollary 3.8.} The equality $\mathrm{Mult}\,A=\mathrm{Mult_0}\,A$ is true if and only if for any $\psi\in\mathrm{Mult}\,A$, we have  $\psi=1$ for any edge $(x,\,y)\in T$ such that $\psi(e_{xy})=e_{xy}$.

\textbf{Corollary 3.9 \cite[Theorem 5]{BruFS15}.} A multiplicative automorphism $\psi$ is inner if and only if for any two elements $x,\,y\in X$, weights of any two semi-paths from $x$ to $y$ coincide.

\textbf{Proof.} If a multiplicative inner automorphism $\psi$ is defined by the element $v=(v_x)$, then the weight of any semi-path from $x$ to $y$ is equal to $v_x^{-1}v_y$. This has already been noted above.

For the automorphism $\psi$, we assume that weight conditions for semi-paths holds. Similar to Theorem 3.6,  we write $\psi=\psi_1\psi_0$. If $(x,\,y)$ is an edge in $T$, then $\psi(e_{xy})=\psi_0(e_{xy})$. Now let $(x,\,y)$ be an edge in $X$ not contained in $T$. Then $\psi(e_{xy})=c_{xy}e_{xy}$ and $w(x,\,y)=c_{xy}$. If $P$ is a semi-path in $T$ from $x$ to $y$, then $w(P)=v_x^{-1}v_y$, where the invertible element $v=(v_x)$ defines an inner automorphism $\psi_0$. It follows from the assumption that $c_{xy}=v_x^{-1}v_y$. Therefore, $\psi$ acts on all elements $e_{xy}$ as $\psi_0$; therefore, $\psi=\psi_0$.~$\square$

Let $C$ be some simple cycle of the graph $X$ and let $\psi$ be some multiplicative automorphism. We define the weight $w(C)$ of the cycle $C$ with respect to $\psi$ as the weight of one of two semi-paths from any vertex of the cycle to itself. When selecting the opposite route direction of the cycle, the weight is equal to $w(C)^{-1}$. Therefore, the weight of the cycle $C$ is one of the elements $w(C)$ and $w(C)^{-1}$.

Let $C_1,\ldots,C_k$ be the fundamental system of cycles associated with  spanning tree $T$ of the graph $X$, where $k=\lambda(X)$.

\textbf{Corollary 3.10.} A multiplicative automorphism $\psi$ is inner if and only if $w(C_1)=\ldots=w(C_k)=1$.

\textbf{Proof.} The necessity of the condition follows from the first paragraph of the proof of Corollary 3.9.

Now we assume that $w(C_1)=\ldots=w(C_k)=1$. We write the multiplicative automorphism $\psi$ in the form $\psi=\psi_1\psi_0$ (see the proof of Corollary 3.9). Let an edge $(x,\,y)$, $x<y$, be not contained in $T$. Then $(x,\,y)\in C_i$ for some $i$. If $P$ is a semi-path in $C_i$ from $x$ to $y$, then we have equalities
$$
1=w(C_i)=w(P)w(x,\,y)^{-1},\quad w(x,\,y)=w(P).
$$
By repeating the argument from the proof of Corollary 3.9, we obtain that $c_{xy}=v_x^{-1}v_y$ and $\psi=\psi_0$.~$\square$

In the book \cite{SpiO97}, multiplicative and fractional automorphisms are defined with the use of some special functions in $I(X,R)$, where $R$ is a commutative ring. These definitions can be also extended to the case of an arbitrary ring $R$.

A function $m\in I(X,R)$ is said to be \textsf{multiplicative} if for every triple of elements $x,\,z,\,y$ in $X$ with $x\leqslant z\leqslant y$, we have
$$
m(x,\,y)\in C(U(R)),\quad m(x,\,y)=m(x,\,z)\cdot m(z,\,y).
$$
Let $m$ be some multiplicative function in $I(X,R)$. Then the mapping
$$
\psi_m\colon A\to A,\; \psi_m(f)=m*f,\; f\in A,
$$
is an automorphism of the algebra $A$ (here $*$ is the Hadamard product of functions, i.e., the pointwise product). It is clear that $\psi_m\in\mathrm{Mult}\,A$. Conversely,  with the use of some multiplicative function, every automorphism in $\mathrm{Mult}\,A$ can be obtained this way.

Let $q\colon X\to C(U(R))$ be some mapping. Then a function $m_q\in A$ such that $m_q(x,\,y)=q^{-1}(x)\cdot q(y)$ for $x\leqslant y$ and $m_q(x,\,y)=0$ in remaining cases, is multiplicative. In addition, the corresponding automorphism $\psi_{m_q}$ is fractional. All fractional automorphisms are obtained in such a way.

\renewcommand{\refname}

\end{document}